\documentclass{amsart}
\usepackage{amsmath}
\usepackage{amsfonts}
\usepackage{amssymb}
\usepackage{graphicx}%
\usepackage{mathrsfs}

\usepackage{times}
\usepackage{bbold}
\usepackage{stmaryrd}
\usepackage{array}
\usepackage{dsfont}
\usepackage[T1]{fontenc}

\newtheorem{theorem}{Theorem}

\newtheorem*{GGR+theorem}{GGR+ Theorem}

\newtheorem{corollary}[theorem]{Corollary}

\newtheorem{definition}[theorem]{Definition}
\newtheorem{example}[theorem]{Example}

\newtheorem{lemma}[theorem]{Lemma}

\newtheorem{proposition}[theorem]{Proposition}

\numberwithin{theorem}{section}
\newcommand{\charf}[1]{\raisebox{\depth}{$\chi$}{\ensuremath{_{#1}}}}

\begin{document}
\title[A generalization of the GGR conjecture]{A generalization of the GGR conjecture}

\author{S. Catoiu}
\address{Department of Mathematics, DePaul University\\Chicago, IL 60614}
\email{scatoiu@depaul.edu}
\author{H. Fejzi\'{c}}
\address{Department of Mathematics, California State University, San Bernardino, CA 92407}
\email{hfejzic@csusb.edu}
\thanks{This paper is in final form and no
version of it will be submitted for publication elsewhere.}
\date{September 14, 2022}

\begin{abstract}
For each positive integer $n$, function $f$, and point $c$, the GGR Theorem states that $f$ is $n$ times Peano differentiable at $c$ if and only if $f$ is $n-1$ times Peano differentiable at $c$ and the following $n$-th generalized Riemann~derivatives of $f$ at $c$ exist:
\[
\lim_{h\rightarrow 0}\frac 1{h^{n}}\sum_{i=0}^n(-1)^i\binom{n}{i}f(c+(n-i-k)h),
\]
for $k=0,\ldots,n-1$. The theorem has been recently proved in \cite{AC2} and has been a conjecture by Ghinchev, Guerragio, and Rocca since 1998.
We provide a new proof of this theorem, based on a generalization of it that produces numerous new sets~of $n$-th Riemann smoothness conditions that can play the role of the above set in the GGR~Theorem.
\end{abstract}
\maketitle

\noindent
Given two finite sequences of real numbers, $\{a_i\}_{i=0}^{m}$ and $\{b_i\}_{i=0}^{m}$, such that the $a_i$ are non-zero and the $b_i$ are distinct, we define
\[
D(h)(f)=\sum_{i=0}^ma_if(c+b_ih),
\]
where $f$ is a real valued function and $c$ is a fixed real number. For a fixed $h$, the map $f\mapsto D(h)(f)$ is a linear operator on the vector space of all real valued functions, while for a fixed $f$, the map $h\mapsto D(h)(f)$ is a function of $h$. The expression $D(h)(f)$ is called a \emph{difference} of $f$ and $h$ centered at $c$, the $a_i$ are its \emph{coefficients}, and the $b_i$ are its \emph{nodes}. For simplicity, we write $D(h)$ to denote $D(h)(f)$ whenever there is no confusion as to what the function $f$ is. We say that the difference $D(h)$ has \emph{order} $n$ if $D(1)(x^k)=0$ for~$k=0,1,\ldots,n-1$ and $D(1)(x^n)\neq 0$. The definition of the order implies that this does not depend on $c$. Moreover, if $D(h)$ has order $n$ then it must have at least $n+1$ nodes and, for every set $B=\{b_i\}_{i=0}^n$, there is a difference of order $n$ with nodes from $B$. It is easy to verify that if $D(h)$ has order $n$, then $D(h)(p)\equiv 0$ for every polynomial $p$ of degree less than~$n$, the value of $D(1)(x^n)$ is independent of $c$, and $D(h)(x^n)=D(1)(x^n)h^n$.

By taking $c=0$, one has $D(1)(x^k)=\sum_{i=0}^ma_i(b_i)^k$, for $k=0,1,\ldots ,n$, so that the difference~$D(h)$ has order $n$ if and only if it satisfies the Vandermonde system of linear equations $\sum_{i=0}^ma_i(b_i)^k=\delta_{k,n}\cdot C$, for~$k=0,1,\ldots ,n$, where $C$ is a non-zero constant. When $C=n!$, the $n$-th difference~$D(h)$ is called an $n$-th \emph{generalized Riemann difference}. In this case, the above linear system in unknowns $a_i$ is consistent when $m\geq n$ and has a unique solution when $m=n$, in which case the $n$-th generalized Riemann difference is called \emph{exact}. Examples of exact $n$-th generalized Riemann differences include the $n$-th Riemann difference $\Delta_n(h)$ and the $n$-th symmetric Riemann difference $\Delta_n^s(h)$, defined by
\[
\text{{\small $\Delta_n(h)=\sum_{i=0}^n(-1)^i\tbinom nif(x+(n-i)h)$} and {\small $\Delta_n^s(h)=\sum_{i=0}^n(-1)^i\tbinom nif(x+(\tfrac n2-i)h)$,}}
\]
the differences in Theorem~\ref{T02} below, and numerous other differences we use in this article.

If $D(h)$ is a difference of order at least $n$, then the limit $\lim_{h\rightarrow 0}D(h)/h^n$ is an $n$-th \emph{smoothness} condition for $f$ at~$c$. When $D(h)$ is an $n$-th generalized Riemann difference, then the $n$-th smoothness condition associated to it is called an $n$-th \emph{generalized Riemann derivative} of $f$ at $c$. When $D(h)$ is just an $n$-th difference, a unique scalar multiple of it is an $n$-th generalized Riemann difference: simply take the scalar to be $n!/C$. Consequently, the $n$-th smoothness condition associated to an $n$-th difference is a non-zero scalar multiple of, hence equivalent to, an $n$-th generalized Riemann derivative. We retain that smoothness is more general than generalized Riemann differentiation.

The generalized Riemann derivatives were introduced by Denjoy in \cite{D} in 1935, generalizing the Riemann derivatives $\mathscr{R}_nf(c)=\lim_{h\rightarrow 0}\Delta_n(h)(f)/h^n$ and the symmetric Riemann derivatives $\mathscr{R}_n^sf(c)=\lim_{h\rightarrow 0}\Delta_n^s(h)(f)/h^n$, invented by Riemann in the mid~1800s; see \cite{R}. Smoothness of order $n$ was introduced by Marcinkiewicz and Zygmund in \cite{MZ} in 1936, and has been recently investigated in \cite{AC2}.

\medskip
The goal of this paper is to use generalized Riemann differences in order to establish sufficient criteria for approximating functions by polynomials. We say that a function $f$ is approximated to order $n$ by a polynomial $p$ near a point $c$, if the error $f(c+h)-p(c+h)$ is of order~$o(h^n)$, that is, if $\lim_{h\rightarrow 0}\frac {f(c+h)-p(c+h)}{h^n}=0$. When this approximation is possible, then there is a unique approximating polynomial of degree less than or equal to~$n$. For~$n$ times differentiable functions at $c$, this approximation is possible and the corresponding result is the well-known Taylor's Theorem. But as it was first pointed by Peano, functions that can be approximated by polynomials do not have to be $n$ times differentiable at $c$. In his honor, we say that $f$ is $n$ times Peano differentiable at $c$, if there is a polynomial $p$ such that~$f(c+h)-p(c+h)=o(h^n)$.

\medskip
Our main result is the following theorem, whose proof is given in Section~\ref{S3}:
\begin{theorem}[The Generalized GGR Theorem]\label{T01}
Let $n\geq 2$ and let $f$ be an $n-1$ times Peano differentiable function at $c$. The following are sufficient conditions for $f$ to be $n$ times Peano differentiable at $c$:
\begin{itemize}
\item For $n$ odd, $\lim_{h\rightarrow 0}D_k(h)(f)/h^n$ exists at $c$, for all $k$ with $\tfrac {n+1}2\leq k\leq n-1$, where $D_k(h)$ is a difference of order at least $n$, whose one node is $k$ and the rest belong to the set $\{-\tfrac {n-1}2,-\tfrac {n-1}2+1,\ldots ,-1,0,1,\ldots ,k-2,k-1\}$.
\item For $n$ even, $\lim_{h\rightarrow 0}D_k(h)(f)/h^n$ exists at $c$, for all $k$ with $\tfrac {n}2\leq k\leq n$, where $D_k(h)$ is a difference of order at least $n$, whose one node is $k$ and the rest belong to the set $\{-\tfrac {n}2,-\tfrac {n}2+1,\ldots ,-1,0,1,\ldots ,k-2,k-1\}$.
\end{itemize}
\end{theorem}

Theorem~\ref{T01} generalizes the GGR Conjecture, which provides a sufficient condition for a function $f$ to be $n$ times Peano differentiable at $c$, using a concrete set of generalized Riemann derivatives instead of a set of smoothness conditions. The conjecture states that the existence of the $\tfrac{n(n+1)}2$ exact generalized Riemann derivatives of $f$ at $c$, of order $k$, for~$k=1,\ldots, n$, and with nodes $\{-j,-j+1,\ldots,-j+k\}$, for $j=0,\ldots,k-2$, is sufficient for $f$ to be $n$ times Peano differentiable at $c$. Using a basic inductive argument, the GGR Conjecture can be expressed in the following equivalent form:

\begin{theorem}[GGR Conjecture]\label{T02}
Let $n\geq 2$ and let $f$ be $n-1$ times Peano differentiable at $c$. If all generalized Riemann derivatives
\[
\text{{\small $\lim_{h\rightarrow 0}\frac {\sum_{j=0}^n(-1)^j\binom {n}{j}f(c+(k-j)h)}{h^n}$}, for $k=1,2,\ldots ,n$, exist}
\]
then $f$ is $n$ times Peano differentiable at $c$.
\end{theorem}

All differences $D_k(h)(f)=\sum_{j=0}^n(-1)^j\binom {n}{j}f(c+(k-j)h)$, for $k=1,2,\ldots ,n$, have order $n$ and, regardless of the parity of $n$, there are more than enough of them to satisfy the hypothesis in Theorem~\ref{T01}, so that the GGR Conjecture is a special case of this theorem. The original GGR Conjecture was proved by Ginchev, Guerragio, and Rocca by hand for $n=1,2,3,4$ in \cite{GGR} and with the help of a computer they proved it for $n=5,6,7,8$ in~\cite{GR}. The conjecture has been recently proved for general $n$ in \cite{AC2} and is now a theorem.

The above exact $n$th generalized Riemann difference, $D_{n-k}$ for $k=0,1,\ldots ,n-1$, was denoted in \cite{AC2} as $\Delta_{n,k}$ and called the \emph{$k$-th backward shift} of the $n$-th Riemann difference $\Delta_n=\Delta_{n,0}$. In this way, the $n$-th symmetric Riemann difference is $\Delta_n^s=\Delta_{n,\frac n2}$.

Our main motivation for generalizing the GGR conjecture in the way we do it in Theorem~\ref{T01} comes from a particular case of a result from \cite{ACCH}, providing two sufficient smoothness conditions for a function at a point that make it differentiable at the point.  We conveniently restate this here in an equivalent form as the following proposition:

\begin{proposition}[\cite{ACCH}, Corollary 3.2]\label{P03}
Given a function $f$ and a point $c$, if both limits
\[
\text{{\small $\displaystyle \lim_{h\rightarrow 0}\frac {f(c+h)-f(c-h)}h$} and {\small $\displaystyle \lim_{h\rightarrow 0}\frac {f(c+h)-2f(c)+f(c-h)}h$} exist,}
\]
then $f$ is differentiable at $c$.
\end{proposition}

Both limits in Proposition~\ref{P03} are smoothness conditions of order 1. The first difference has order~1 and~2 nodes; its smoothness condition is twice the symmetric derivative of~$f$ at~$c$. The second difference has order~2.

When $n=1$, Theorem~\ref{T02} is a tautology that does not even require the hypothesis that~$f$ is $n-1$ times Peano differentiable at $c$; while Theorem~\ref{T01} does not make sense. Proposition~\ref{P03} is almost an extension of Theorem~\ref{T01} for $n=1$; it points in the direction of the result for higher $n$. A non-trivial generalization of the $n=1$ case in Theorem~\ref{T02}, the problem of finding all exact first order generalized Riemann derivatives that imply the first (Peano) derivative for all continuous functions at $c$, is proved in \cite{C}.

\medskip
Later on in the introduction we will explain why the condition that $f$ is $n-1$ times Peano differentiable at $c$ in Theorem~\ref{T01} is necessary. Before that, we provide a few examples for small $n$, to help with understanding the result of the theorem. In addition, Examples~\ref{E02}, \ref{E04}, \ref{E05}, and \ref{E07} justify the need for more than one smoothness condition in the theorem, a need not highlighted in any of the above mentioned articles.

In the case of Proposition~\ref{P03}, consider the following two non-differentiable functions at $c$: $f(x)=|x-c|$, for which the first smoothness condition is satisfied while the second is not; and $g(x)=|x-c|\cdot \charf{[c,\infty)}(x)$, where $\charf{[c,\infty)}$ is the characteristic function of the interval $[c,\infty )$, for which the second smoothness condition in the proposition is satisfied and the first is not. The two examples highlight the fact that neither of the two smoothness conditions alone implies the differentiability of a function at $c$, but together they do.

\medskip
Here are the announced examples:

\begin{example}\label{E02}
{\rm
Let $n=2$ and let $f$ be differentiable at $c$. By Theorem~\ref{T01}, if both limits
\[
\text{{\small $\displaystyle \lim_{h\rightarrow 0} \frac {f(c+h)-2f(c)+f(c-h)}{h^2}$} and {\small $\displaystyle \lim_{h\rightarrow 0} \frac {f(c+2h)-3f(c+h)+3f(c)-f(c-h)}{h^2}$}}
\]
exist, then $f$ is twice Peano differentiable at $c$.
}
\end{example}
The above example has $1\leq k\leq 2$, so that $D_1(h)$ has to be a difference of order~2 with nodes $\{-1,0,1\}$, and we picked the exact difference $\Delta_2^s$; while for $D_2(h)$ we picked the third difference $\Delta_{3,1}$ out of a wider range of differences of orders at least 2 and with nodes $\{-1,0,1,2\}$. That the first limit condition alone is not enough to assure that $f$ is twice Peano differentiable at $c=0$ is easily seen by analyzing the function $f(x)=x|x|$, which is differentiable at 0 and has $D_1(h)(f)\equiv 0$, while $f$ is not twice Peano differentiable at 0.

\begin{example}\label{E03}
{\rm
Let $n=3$ and let $f$ be twice Peano differentiable at $c$. By Theorem~\ref{T01}, if
\[
\text{{\small $\displaystyle \lim_{h\rightarrow 0} \frac {f(c+2h)-3f(c+h)+3f(c)-f(c-h)}{h^3}$} exists,}
\]
then $f$ is three times Peano differentiable at $c$.
}
\end{example}
Given that Example~\ref{E02} needed two limits to assure twice Peano differentiation, Example~\ref{E03} requiring a single limit to guarantee third Peano differentiation comes as a surprise. Indeed, the condition $\tfrac {3+1}2\leq k\leq 3-1$ in Theorem~\ref{T01} allows only one choice for $k$. Other instances where a single such limit implies three times Peano differentiation for all twice Peano differentiable functions $f$ at $c$ are provided in \cite{ACF,ACF1}.

\begin{example}\label{E04}
{\rm
Let $n=4$ and let $f$ be three times Peano differentiable at $c$. If all of the following limits

{\small \[
\begin{aligned}
&\lim_{h\rightarrow 0} [{f(c+2h)-4f(c+h)+6f(c)-4f(c-h)+f(c-2h)}]/{h^4},\\
&\lim_{h\rightarrow 0} [{f(c+3h)-4f(c+2h)+6f(c+h)-4f(c)+f(c-h)}]/{h^4},\\
&\lim_{h\rightarrow 0} [{f(c+4h)-4f(c+3h)+6f(c+2h)-4f(c+h)+f(c)}]/{h^4},
\end{aligned}
\]}
exist, then $f$ is four times Peano differentiable at $c$.
}
\end{example}

In Example~\ref{E04}, by Theorem~\ref{T01}, $n=4$ implies $2\leq k\leq 4$, so that $D_2(h)$ has to be a difference of order 4 with nodes $\{-2,-1,0,1,2\}$, hence a nonzero scalar multiple of~$\Delta_{4,2}$, and we picked the scalar to be 1; while for $D_3(h)$ and $D_4(h)$ we picked the fourth differences $\Delta_{4,1}$ and $\Delta_4$ out of two wider classes of differences of orders at least 4 whose respective sets of nodes are included in $\{-2,-1,\ldots,3\}$ and $\{-2,-1,\ldots ,4\}$.

The following example shows that only the first two limit conditions in Example~\ref{E04}, those corresponding to the differences $D_2(h)$ and $D_3(h)$, are not enough to guarantee four times Peano differentiability for all three times Peano differentiable functions at $c$.

\begin{example}\label{E05}
{\rm Let $s$, with $3<s<4$, be a real number. For $h\geq 0$, define $f(c+h)=(-1)^{m+n}h^s$ if $h=2^m3^n$ where $m$ and $n$ are integers, and 0 otherwise. For $h<0$ we define $f(c+h)=-f(c-h)$. Clearly $f$ is three times Peano differentiable at $c$, but not four times.

Since $f$ is odd relative to $c$ and $f(c)=0$, $D_2(h)=0$, while $D_3(h)=f(c+3h)-4f(c+2h)+5f(c+h)=(-1)^{m+n}h^s(-3^s+4\times 2^s+5)$, for $h>0$. Since $-3^3+4\times 2^3+5=10>0$ and $-3^4+4\times 2^4+5=-12<0$, if we pick $s$ with $3<s<4$ such that $-3^s+4\times 2^s+5=0$, then $D_3(h)\equiv 0$ for $h>0$, and the property that $D_3(-h)(f)=-D_3(h)(f)$ makes $D_3(h)\equiv 0$ for all $h$.
}
\end{example}

\begin{example}\label{E06}
{\rm
Let $n=5$ and let $f$ be four times Peano differentiable at $c$. If both limits
{\small \[
\begin{aligned}
&\lim_{h\rightarrow 0} [{f(c+3h)-5f(c+2h)+10f(c+h)-10f(c)+5f(c-h)-f(x-2h)}]/{h^5}\text{ and}\\
&\lim_{h\rightarrow 0} [{f(c+4h)-15f(c+2h)+40f(c+h)-45f(c)+24f(c-h)-5f(c-2h)}]/{h^5}
\end{aligned}
\]}
exist, then $f$ is five times Peano differentiable at $c$.
}
\end{example}

Theorem~\ref{T01} is applied in Example~\ref{E06} for $3\leq k\leq 4$, so that $D_3(h)$ has to be a difference of order 5 with nodes $\{-2,-1,0,1,2,3\}$, hence a non-zero scalar multiple of~$\Delta_{5,2}(h)$, and we considered the case when the scalar is 1; for $D_4(h)$ we picked a difference of order 5 out of a wider range of possibilities of differences of orders at least 5 and with set of nodes included in $\{-2,-1,\ldots ,4\}$.

The next example shows that the first limit condition in Example~\ref{E06} is not enough to imply five times Peano differentiability for all four times Peano differentiable functions~at~$c$.

\begin{example}\label{E07}
{\rm Let $s$, with $4<s<5$, be a real number. For $h\geq 0$, define $f(c+h)=(-1)^{m+n}h^s$ if $h=2^m3^n$, where $m$ and $n$ are integers, and 0 otherwise. For $h<0$ we define $f(c+h)=f(c-h)$. Clearly $f$ is four times Peano differentiable at $c$, but not five times Peano differentiable at $c$.

Since $f$ is even relative to $c$ and $f(c)=0$, for $h>0$, $D_3(h)=f(c+3h)-6f(c+2h)+15f(c+h)=(-1)^{m+n}h^sp(s)$, where $p(s)=-3^s+6\times 2^s+15$. And since $p(4)=30>0$ and $p(5)=-36<0$, there is an $s$ between 4 and 5 such that $p(s)=0$, so that $D_3(h)\equiv 0$ for $h>0$. This extends to $D_3(h)\equiv 0$ for all $h$, due to the hypothesis that~$f$ is even relative to $c$.
}
\end{example}

As promised, we will show that the condition that $f$ is $n-1$ times Peano differentiable at~$c$ in the statement of Theorem~\ref{T01} is necessary.
This will follow from Theorem~\ref{T08} below, which also highlights the following two consequences of Theorem~\ref{T01}.
\begin{itemize}
\item It is easy to see that the difference $D_k(h)$ corresponding to the lowest $k$ in Theorem~\ref{T01} has $n+1$ nodes and order at least $n$, so it must have order $n$.
\item In general, if a finite collection $\{D_{\alpha }(h)\}$ of differences of orders at least $n$ has the property that
$\lim_{h\rightarrow 0}D_{\alpha }(h)(f)/h^n$ exists for all $\alpha $ implies that $f$ is $n$ times Peano differentiable at $c$, for all $n-1$ times Peano differentiable functions $f$ at $c$, then at least one $D_{\alpha }(h)$ must have order $n$.
\end{itemize}

\begin{theorem}\label{T08}
Let $n$ be an integer, at least 2, and let $\{D_{\alpha }(h)\}$ be a finite collection of differences of orders at least $n$. Then:
\begin{enumerate}
\item[(i)] There is an $n-2$ times Peano differentiable function~$f$ at $c$, so that $D_{\alpha }(h)(f)\equiv 0$ for all $\alpha$, but $f$ is not $n$ times Peano differentiable at $c$.
\item[(ii)] If all orders are greater than $n$, then there is an $n-1$ times Peano differentiable function $f$ at $c$, such that $D_{\alpha }(h)(f)\equiv 0$ for all $\alpha$, but $f$ is not $n$ times Peano differentiable~at~$c$.
\end{enumerate}
\end{theorem}

\begin{proof}
Let $B=\{b_j\}_{j=1}^m$ denote the set of all nodes of all the $D_{\alpha}(h)$. The set $G=\{\prod_{j=1,b_j\neq 0}^mb_j^{k_j}\mid k_j\in\mathbb{Z}\}$ has the property that if $h\in G$ then $b_jh\in G$ for all $0\neq b_j\in B$, and if $h\notin G$ then $b_jh\notin G$ for all $b_j\in B$. Let $f,g:\mathbb{R}\rightarrow \mathbb{R}$, defined~by
\[
f(c+h)=\begin{cases} h^{n-1}&\text{if }h\in G\\ 0 &\text{otherwise}\end{cases}\quad\text{ and }\quad
g(c+h)=\begin{cases} h^{n}&\text{if }h\in G\\ 0 &\text{otherwise}\end{cases}.
\]
Clearly, $B$ finite makes $G$ countable, so that both $f$ and $g$ are measurable. It is also clear that $f$ is $n-2$ and $g$ is $n-1$ times Peano differentiable at $c$, but neither $f$ nor $g$ are $n$ times Peano differentiable at $c$.

(i) Since each difference of order greater than $n-1$ vanishes on every polynomial degree up to $n-1$, $D_{\alpha }(h)(f)\equiv 0$ for all $\alpha$. (ii) If all orders are greater than $n$, then a similar argument for $n$ in place of $n-1$ makes $D_{\alpha }(h)(g)\equiv 0$ for all $\alpha$.
\end{proof}

The proof of the GGR Theorem given in \cite{AC2} has a part based on the theory of symmetric Peano and symmetric generalized Riemann derivatives, developed in that article, and a part based on a highly non-trivial combinatorial algorithm. The proof of the Generalized GGR Theorem, Theorem~\ref{T01}, is based entirely on analysis, by extending the notion of a generalized Riemann differentiation to the notion of a generalized Riemann smoothness. It uses properties of symmetric differences, which we review next.

\subsection*{Symmetric differences} A difference $D(h)(f)$ of a function $f(x)$ is an \emph{even or odd difference}, if $D(-h)=\pm D(h)$, and is a \emph{symmetric difference}, if $D(-h)=(-1)^n D(h)$, where $n$ is the order of $D$. For example, the $n$-th symmetric Riemann difference $\Delta_n^s(h)$ is a symmetric difference. The following are properties of symmetric differences that will be used throughout the paper:
\begin{itemize}
\item The set of nodes of a symmetric difference is symmetric relative to the origin.
\item Odd differences do not allow 0 as a node.
\item The order of a symmetric difference has the same parity as the difference.
\item Each exact difference whose set of nodes is symmetric relative to the origin must be a symmetric difference.
\item A linear combination of symmetric differences of the same parity is a symmetric difference of the same parity.
\end{itemize}
Most of these properties can be checked directly, from the definition of a symmetric difference; others have a bit more involved proofs. For more on symmetric differences, see~\cite{AC2,ACCh}.

\medskip
The exact even symmetric difference $S_k$ with nodes $0,\pm1,\pm2,\ldots,\pm k$ is the difference \[S_k=\Delta_{2k}^s=\Delta_{2k,k}.\] The exact odd symmetric difference $T_k$  with nodes $\pm1,\pm2,\ldots,\pm k$ is not the difference $\Delta_{2k-1}^s=\Delta_{2k-1,(2k-1)/2}$, whose nodes are half integers. Its actual expression is given by
\[
T_k=\left(\Delta_{2k-1,k-1}+\Delta_{2k-1,k}\right) /2.
\]

\subsection*{A few details of the proof} 
The proof of the Generalized GGR Theorem relies on two fundamental lemmas that are the subject of Sections~\ref{S1} and~\ref{S2}.

The first fundamental lemma, Theorem~\ref{th2}, is needed in the proof of the second fundamental lemma. It says that if $R(h)=o(h^s)$ and $D(h):=R(2h)-2^sR(h)=o(h^n)$, for $n>s\geq 0$, then $R(h)=o(h^n)$. It is worded in terms of functions, but it can also be worded in terms of differences of the same function at $x=0$ and $h$.

The second fundamental lemma is a result on symmetric differences of $f$ at $0$, where~$f$ satisfies $f(h)=o(h^{n-1})$. Since it has many ingredients that are parity dependent, this lemma is separated into two theorems with similar proofs: Theorem~\ref{th3} for even differences, and Theorem~\ref{th4} for odd differences.

The statement of the second fundamental lemma is closely related to the statement and especially the original proof of the GGR Theorem in \cite{AC2}, that uses symmetric differences. It says that when the differences $D_k,D_{k+1},\ldots, D_{n-1},(D_n)$ of $f$ at 0 are symmetric (of a more general kind when compared to the ones in the original proof of the GGR Theorem), have the same parity, and have an extra condition on their orders, then this sequence can be completed to a sequence with the same properties all the way down to either $k=2$ or~$1$, depending on the common parity of the differences. In particular, this leads to the existence of one of the following two limits:
\[
\begin{aligned}
&\lim_{h \to 0}{[f(h)+f(-h)]}/{h^n}\text{,\qquad when all the $D_j$ are even;}\\
&\lim_{h \to 0}{[f(h)-f(-h)]}/{h^n}\text{,\qquad when all the $D_j$ are odd.}
\end{aligned}
\]

The last step between the second fundamental lemma and the proof of the Generalized GGR Theorem is Proposition~\ref{cr1}, dealing with the case when there are two such sequences of differences $D_j$ of $f$ at $x$, one of even differences and one of odd differences, all under the assumption that $f$ is $n-1$ times Peano differentiable at $x$. In this case, 
by subtracting from $f(x)$ its $(n-1)$-st Taylor polynomial and then by shifting the variable, without loss of generality, we may assume that $x=0$, $f(0)=0$, and $f(h)=o(h^{n-1})$. Then the second fundamental lemma makes both of the above limits exist and so does their average, $\lim_{h \to 0}{f(h)}/{h^n}$. This in turn is equivalent to the $n$-th Peano differentiability of $f$ at $x$.

The leftover proof of the Generalized GGR Theorem is then about using the hypothesis to construct the two sequences $\{D_j\}$ in Proposition~\ref{cr1}, one for each parity differences.


\subsection*{The context} 
The solution to the GGR Conjecture came in the context of a few recent developments in the theory of generalized derivatives. Among these are the solution in \cite{ACCs}, by Ash, Catoiu, and Cs\"ornyei, to the nearly a century old problem of finding the equivalences between Peano and generalized Riemann derivatives, and the solution in \cite{ACCh,ACCH}, by Ash, Catoiu, and Chin, to the half a century old problem of finding the equivalences between any two generalized Riemann derivatives. 
These developments highlighted the use of new techniques on generalized derivatives that are usually proper to other areas of mathematics, such as the theory of infinite linear systems, the use of group algebras, grading, or generalized polynomial algebras, recursive set theory and combinatorial methods.

The problem of finding the equivalences between the Peano and generalized Riemann derivatives, or the above mentioned century old problem, was initiated in 1927 by Kintchine in \cite{Ki}, who proved that the symmetric derivative is equivalent to the first Peano derivative, for all functions $f$ at almost everywhere points $c$ on a measurable set. This was extended to order $n$ symmetric Riemann and Peano derivatives in~1936 by Marcinkiewicz and Zygmund in \cite{MZ}, and then to order $n$ generalized Riemann and Peano derivatives in 1967 by J. Marshall Ash in \cite{As}. The original proofs in \cite{As} and \cite{MZ} assumed that certain sets were measurable. This flaw was corrected by Fejzi\'c and Weil in \cite{FW1}.

The generalized Riemann derivatives were shown to satisfy properties similar to those for ordinary derivatives, such as monotonicity \cite{HL1,T,W}, convexity \cite{GGR1,HL,MM}, or the mean value theorem \cite{AJ,FFR}. They have many applications in the theory of trigonometric series \cite{SZ,Z} and numerical analysis \cite{AJJ,L,Sa}. Quantum Riemann derivatives are studied in \cite{AC,ACR}, and multidimensional Riemann derivatives are a part of \cite{AC1}. 
For more on generalized Riemann differentiation see \cite{RAA1,RAA2} and the survey article~\cite{As2}~by~Ash.

The study of Peano derivatives also has a long and rich history. See for example the survey article \cite{EW} by Evans and Weil.
The Peano derivatives were invented by Peano in~\cite{P} and then developed by de la Vall\'ee Poussin in \cite{dlVP}.
More recent developments in this subject can be found in \cite{F,F1,FR,FW}.

\section{The first fundamental lemma}\label{S1}

The following theorem is crucial in the proofs of the main results in Section~\ref{S2}.

\begin{theorem}\label{th2} Let $n>s$ be non-negative integers, and let $D(h)=R(2h)-2^sR(h)$, where $R$ is a function of~$h$.
If $R(h)=o(h^s)$ and $D(h)=o(h^n)$ then~$R(h)=o(h^n)$.
\end{theorem}
 
\begin{proof}
By the hypothesis $D(h)=o(h^n)$, there are $\epsilon ,\delta >0$ such that $|D(h)|<\epsilon |h|^n$, whenever $0<|h|<\delta$. Fix $h$. Substituting  $t=\frac{h}{2^k}$, for $k=1,2,\ldots,\ell $, into the hypothesis
$$D(t)=R(2t)-2^sR(t),$$ and observing that $|t|=|\frac{h}{2^k}|<\delta$ makes $|D(t)|<\epsilon |t|^n$, we obtain a system of inequalities,
$$
 \begin{aligned}
 | R(h)-2^sR\left(\tfrac{h}{2}\right)|<\epsilon \tfrac{ |h|^n}{2^n},&\quad
 |R\left(\tfrac{h}{2}\right)-2^sR\left(\tfrac{h}{2^2}\right)|<\epsilon\tfrac{ |h|^n}{2^{2n}},\ldots \\
  \ldots ,\;| R\left(\tfrac{h}{2^{\ell -1}}\right)-&2^sR\left(\tfrac{h}{2^{\ell }}\right)| <\epsilon \tfrac{ |h|^n}{2^{\ell n}}.
 \end{aligned}
 $$
Multiply the second inequality by $2^s$, the third by $2^{2s}$, and so on, the last one by~$2^{(\ell -1)s}$, to deduce
 
 \begin{equation}\label{eq0}
 \begin{aligned}
 | R(h)-2^sR\left(\tfrac{h}{2}\right)|<\epsilon \tfrac{ |h|^n}{2^n},&\quad
 |2^sR\left(\tfrac{h}{2}\right)-2^{2s}R\left(\tfrac{h}{2^2}\right)|<\epsilon\tfrac{ |h|^n2^s}{2^{2n}},\ldots \\
 \ldots , |2^{(\ell -1)s} R\left(\tfrac{h}{2^{\ell -1}}\right)&-2^{\ell s}R\left(\tfrac{h}{2^{\ell }}\right)| <\epsilon \tfrac{ |h|^n2^{(\ell -1)s}}{2^{\ell n}}.
 \end{aligned}
 \end{equation}
 Since
$$
\begin{aligned}
R(h)-2^{\ell s}R\left(\tfrac{h}{2^{\ell }}\right)&=
[R(h)-2^sR\left(\tfrac{h}{q}\right)]+[2^sR\left(\tfrac{h}{2}\right)-2^{2s}R\left(\tfrac{h}{2^2}\right)]+\cdots\\
&\cdots+[2^{(\ell -1)s}R\left(\tfrac{h}{2^{\ell -1}}\right)-2^{\ell s}R\left(\tfrac{h}{2^{\ell }}\right)],
\end{aligned}
 $$ 
by the triangle inequality, the sum of all inequalities in \eqref{eq0} implies the single inequality
 $$| R(h)-2^{\ell s}R\left(\tfrac{h}{2^{\ell }}\right)|<\epsilon\tfrac{ |h|^n}{2^{n}}\sum_{k=0}^{\ell -1} \left(\tfrac{2^s }{2^{n}}\right)^k,$$ and from here, 
\begin{equation}\label{eq1a}
| R(h)|<2^{\ell s}|R\left(\tfrac{h}{2^{\ell }}\right)| +\epsilon\tfrac{ |h|^n}{2^{n}}B,
\end{equation}
where $B$ is the sum of the convergent geometric series $\sum_{k=0}^{\infty} \left(\frac{2^s }{2^{n}}\right)^k$. 
To complete the proof, it suffices to show that the term $2^{\ell s}R\left(\frac{h}{2^{\ell s}}\right)$ in \eqref{eq1a} goes to 0 as $\ell  \to \infty$. Indeed,~we $$\text{rewrite this as }h^s{R\left(\tfrac{h}{2^{\ell }}\right)}/{(\tfrac{h}{2^{\ell }})^s},$$
and the desired assertion follows from the hypothesis that $R(t)=o(t^s)$.
\end{proof}

Note that both the statement and the proof of Theorem~\ref{th2} will remain true if the number~2 in the definition of $D(h)$ is be replaced by any real number $q$ greater than 1.

\section{The second fundamental lemma}\label{S2}

The second fundamental lemma is divided into two parts: one result for even differences, Theorem~\ref{th3}, and another result for odd differences, Theorem~\ref{th4}. The proofs of both theorems require the following lemma on basic properties of the two special symmetric differences defined in the introduction: the exact even generalized Riemann difference $S_k=\Delta_{2k,k}$ with nodes $0,\pm1,\ldots,\pm k$ and order $2k$; and the exact odd generalized Riemann difference $T_k$ with nodes $\pm1,\pm2,\ldots,\pm k$ and order $2k-1$.
Before stating the lemma, we need to make two remarks.

The first remark is about two interpretations of the Vandermonde relations of a generalized Riemann difference, that will be used in both the statement and the proof of Lemma~\ref{pr60}, as well as in the proofs of the two versions of the second fundamental lemma. Suppose $D(h)(f)=\sum_{i}A_if(x+a_ih)$ is such a difference of order $p$. The left side of its~$j$th Vandermonde relation can be interpreted as
\[
\text{$\sum_iA_i(a_i)^j$ is equal to $D(1)(x^j)$ evaluated at $x=0$,}
\]
and the trivial equation $\sum_iA_i(a_i)^j=\tfrac 1{h^j}\sum_iA_i(a_ih)^j$ is interpreted in the original difference as
\[
\text{$D(1)(x^j)=D(h)(x^j)/h^j$, for all $h$, when both sides are evaluated at $x=0$.}
\]

The second remark is an interpretation using Laurent polynomials of the differences that appear in the proof of Lemma~\ref{pr60}. The map that associates to a difference $D(h)=\sum_iA_if(x+a_ih)$ the Laurent polynomial $p(t)=\sum_iA_it^{a_i}$ is a linear isomorphism from the space of all differences with integer nodes to the space of all Laurent polynomials. Here are a few properties of this isomorphism: if $D(h)=\Delta_{n,j}(h)$, then $p(t)=t^{-j}(t-1)^n$; $D(h)$ is an $n$-th generalized Riemann difference, if $p(t)=(t-1)^nq(t)$ and $q(1)=1$, and; $D(h)$ is an~$n$-th difference, if $p(t)=(t-1)^nq(t)$ and $q(1)$ is non-zero.

\medskip
The following lemma uses the above interpretation for $D(h)$ when this is either~$S_k$~or~$T_k$.
\begin{lemma}\label{pr60}
For a positive integer $k$, the differences $S_k$ and~$T_k$ satisfy the following properties:
\begin{enumerate}
\item[(i)\,] $S_k(1)(x^n)\neq 0$, for all even $n$, with $n\geq 2k$;
\item[(ii)] $T_k(1)(x^n)\neq 0$, for all odd $n$, with $n\geq 2k-1$.
\end{enumerate}
\end{lemma}

\begin{proof}
(i) By the above two remarks, the polynomial $p(t)$ associated to the difference $S_k(h)=\Delta_{2k,k}(h)$ is $p(t)=t^{-k}(t-1)^{2k}=(t^{\frac 12}-t^{-\frac 12})^{2k}=\sum_{j=-k}^ka_jt^j$. Consider the function $P(x):=p(e^x)=(e^{\frac x2}-e^{-\frac x2})^{2k}=x^{2k}(c_0+c_2x^2+c_4x^4+\cdots )^{2k}$, for some positive constants $c_j$. It follows that for $n$ even, $n\geq 2k$, the coefficient of $x^n/n!$ in the Taylor series expansion of $P(x)$ is not zero. Since $P(x)=\sum a_je^{jx}$, the same coefficient is~$P^{(n)}(0)=\sum a_j j^n=S_k(1)(x^n)$, and the result is clear.

(ii) In this case, by the same remarks, the polynomial associated to the difference $2T_k(h)=\Delta_{2k-1,k-1}(h)+\Delta_{2k-1,k}(h)$ is $p(t)=t^{-k+1}(t-1)^{2k-1}+t^{-k}(t-1)^{2k-1}=t^{-k+1}(t-1)^{2k-2}(t-t^{-1})=(t^{\frac 12}-t^{-\frac 12})^{2k-2}(t-t^{-1})
=2\sum a_j t^j$. Then $P(x):=p(e^x)=(e^{\frac x2}-e^{-\frac x2})^{2k-2}(e^x-e^{-x})=(c_1x+c_3x^3+c_5x^5+\cdots )^{2k-2}(d_1x+d_3x^3+d_5x^5+\cdots)$ has the $c_i$ and the $d_i$ positive, hence non-zero, for all odd $i$. In particular, for odd $n$, with $n\geq 2k-1$, the coefficient of $x^n/n!$ in the Taylor expansion of $P(x)$ is non-zero. Writing $P(x)=2\sum a_je^{jx}$, this coefficient is $2\sum a_jj^n=2T_k(1)(x^n)$, yielding the result.
\end{proof}
\subsection{The second fundamental lemma for even differences}

\begin{definition}\label{d1}
{\rm Let $n,k$ be positive integers, with $n\geq 3$, and let $m=\lfloor\tfrac{n-1}2\rfloor $. For a function $f$ and point $x$, a set $\{D_k, D_{k+1},\ldots,D_{2m}\}$ of even differences of $f$ at $x$ is called an \emph{$\mathcal{S}_k$-class of even differences}, or an \emph{$\mathcal{S}_k^{{\rm ev}}$-class}, if, for each $r=k,k+1,\ldots,2m$,
\begin{enumerate}
\item[(i)\;] $D_r$ is an even difference of order at least $2k$;
\item[(ii)\,] two nodes of $D_r$ are $\pm r$ and the rest belong to the set $\{0,\pm 1,\ldots,\pm (r-1)\}$;
\item[(iii)] $\lim_{h\rightarrow 0}D_r(h)(f)/h^n$ exists.
\end{enumerate}
For convenience, we also denote $\mathcal{S}_k=\mathcal{S}_k^{\rm ev}=\{D_k, D_{k+1},\ldots,D_{2m}\}$ whenever the above conditions are satisfied.
}
\end{definition}

The following is the second fundamental lemma for even differences.

\begin{theorem}\label{th3}
Suppose that $f(h)=o(h^n)$. Then each $\mathcal{S}_{m+1}$-class of even differences of $f$ at 0 extends to an $\mathcal{S}_1$-class. In particular, the limit $$\lim_{h \to 0}{[f(h)+f(-h)]}/{h^n}\text{ exists.}$$
\end{theorem}

\begin{proof}
Let $k$ be the smallest positive integer at most $m$ such that there is an $\mathcal{S}_{k+1}$-class of even differences. We have to show that $k=0$. If $k>0$, it suffices to show that $\lim_{h\rightarrow 0}S_k(h)/h^n$ exists. Then the set $\{S_k\} \bigcup \mathcal{S}_{k+1}$ will be an $\mathcal{S}_k$-class of even differences, contradicting the minimality of~$k$.
Consider the even difference $$S(h):=S_k(2h)-2^{2k}S_k(h)$$ and claim that $S(h)$ has order at least $2k+2$. To prove this, it suffices to show that $S(1)(x^p)=0$, for $p=0,1,\ldots,2k+1$. Indeed, for $p\leq 2k$, $S(1)(x^p)=S_k(1)((2x)^p)-2^{2k}S_k(1)(x^p)=(2^p-2^{2k})S_k(1)(x^p)=(2^p-2^{2k})\delta_{p,2k}\cdot (2k)!=0$, where the second to last equality follows from the Vandermonde relations for $S_k$. And, for $p=2k+1$, $S(1)(x^p)=0$ is due to $S(h)$ being an even difference.

All nodes of $S(h)$ belong to $\{0,\pm 1,\ldots,\pm k,\ldots ,\pm 2k\}$. By the assumption~$k\leq m$, this is a subset of $\{0,\pm 1,\ldots,\pm k,\ldots ,\pm 2m\}$. Since $\mathcal{S}_{k+1}=\{D_{k+1},\ldots,D_{2m}\}$ and the largest node in $D_r(h)$ is~$r$, there exist $c_{k+1},c_{k+2},\ldots, c_{2m}$ such that all nodes of the difference

$$S(h)-\sum_{ D_r\in \mathcal{S}_{k+1}} c_r D_r(h)$$ belong to the set $\{-2m,-2m+1,\ldots, k-1,k\}$. This together with the fact that this is an even difference makes all nodes belong to $\{0,\pm 1,\ldots,\pm k\}$. And since the order of the difference is at least $2k+2$, that is, more than the number of its base points, this difference has to be the zero difference. Consequently, $S(h)=\sum_{ D_r\in \mathcal{S}_{k+1}} c_r D_r(h)$ and so the~limit  
\begin{equation}\label{eq1}
\lim_{h \to 0}{S(h)(f)}/{h^n}=K \text{ exists.}
\end{equation}
\\
If $n$ is odd, then replacing $h$ with $-h$ in the above limit leads to $-K=K$, or $K=0$.
If $n$ is even, the condition $n>2m\geq 2k$ implies by Lemma~\ref{pr60} that $L:=S_k(1)(x^n)$ is non-zero, which in turn leads to $S(1)(x^n)={S(h)(x^n)}/{h^n}$
{\small
$$
=\frac {S_k(2h)(x^n)-2^{2k}S_k(h)(x^n)}{h^n}=\frac {S_k(1)(x^n)(2h)^n-2^{2k}S_k(1)(x^n)h^n}{h^n}=(2^n-2^{2k})L=M
$$}is non-zero. Consider the function $g$ defined by
$$g(x)=\begin{cases}f(x)-(K/M)x^n & \text{  if $n$ is even,}\\
f(x) & \text{  if $n$ is odd.}
\end{cases}
$$
In both parity cases, equation \eqref{eq1} implies that $S(h)(g)=o(h^n)$.

We need to apply Theorem~\ref{th2}, with $D(h)=S(h)(g)$, $R(h)=S_k(h)(g)$, and $s=2k$. For this we have to check that  $S_k(h)(g)=o(h^{2k})$. Indeed, this follows from the expression
  
\begin{equation}\label{eq2}
S_k(h)(g)=\begin{cases}S_k(h)(f)-(KL/M)h^n & \text{  if $n$ is even,}\\
S_k(h)(f) & \text{  if $n$ is odd,}
\end{cases}
\end{equation}
since the inequality $2k\leq  2m\leq n-1$ together with the hypothesis $f(h)=o(h^{n-1})$ gives $f(h)=o(h^{2k})$, and so $S_k(h)(f)=o(h^{2k})$, and clearly $h^n=o(h^{2k})$.
Theorem~\ref{th2} makes $S_k(h)(g)=o(h^n)$, and so by~\eqref{eq2}, the limit $$\lim_{h \to 0}\frac{S_k(h)(f)}{h^n}=\begin{cases}KL/M & \text{ if $n$ is even,}\\0 & \text{ if $n$ is odd,}\end{cases}$$
exists, as needed.
\end{proof}

\subsection{The second fundamental lemma for odd differences}
The fundamental result in this case operates with the notion of an $\mathcal{S}_k$-class of odd differences, which we define next.

\begin{definition}\label{d2}
{\rm Let $n,k$ be positive integers, with $n\geq 3$, and let $m=\lfloor\tfrac{n}2\rfloor $. For a function~$f$ and point $x$, a set $\{D_k, D_{k+1},\ldots,D_{2m}\}$ of odd differences of $f$ at $x$ is an \emph{$\mathcal{S}_k$-class of odd differences}, or an \emph{$\mathcal{S}_k^{\rm odd}$-class}, if, for each $r=k,k+1,\ldots,2m$,
\begin{enumerate}
\item[(i)\;] $D_r$ is an odd difference of order at least $2k-1$;
\item[(ii)\,] two nodes of $D_r$ are $\pm r$ and the rest belong to the set $\{\pm 1,\ldots,\pm (r-1)\}$;
\item[(iii)] $\lim_{h\rightarrow 0}D_r(h)(f)/h^n$ exists.
\end{enumerate}
Again, for convenience, we also write $\mathcal{S}_k=\mathcal{S}_k^{\rm odd}=\{D_k, D_{k+1},\ldots,D_{2m}\}$ whenever the above conditions are satisfied.
}
\end{definition}

The following is the second fundamental lemma for odd differences.

\begin{theorem}\label{th4}
Suppose that $f(h)=o(h^{n-1})$. Then each $\mathcal{S}_{m+1}$-class of odd differences of $f$ at 0 extends to an $\mathcal{S}_1$-class. In particular, the limit $$\lim_{h \to 0}{[f(h)-f(-h)]}/{h^n}\text{ exists.}$$
\end{theorem}

\begin{proof}
Let $k\leq m$ be the smallest positive integer for which there is an $\mathcal{S}_{k+1}$-class of odd differences, and show that $k=0$. If $k>0$, it is enough to prove that the limit $\lim_{h \to 0}{T_k(h)(f)}/{h^n}$ exists, for then $\{T_k\} \bigcup \mathcal{S}_{k+1}$ is an $\mathcal{S}_k$-class of odd differences, and this would contradict the minimality of $k$. Similar to the proof of Theorem~\ref{th3}, one can show that the odd difference $$T(h)=T_k(2h)-2^{2k-1}T_k(h)$$ has order at least $2k+1$ and its nodes belong to the set $\{\pm 1,\ldots,\pm k,\ldots ,\pm 2k\}$,
hence they belong to $\{\pm 1,\ldots,\pm k,\ldots ,\pm 2m\}$. Moreover, there exist coefficients $c_{k+1},c_{k+2},\ldots,c_{2m}$ such that all nodes of the difference 

$$T(h)-\sum_{ D_r\in \mathcal{S}_{k+1}} c_r D_r(h)$$ belong to the set $\{ -2m,\ldots ,-k,\ldots, -1,1,\ldots ,k\}$, and, by the symmetry of the odd difference, they must belong to $\{\pm 1,\ldots,\pm k\}$. With order at least $2k+1$ and having at most $2k$ nodes, the above difference must be the zero difference. It follows that $T(h)=\sum_{ D_r\in \mathcal{S}_{k+1}} c_r D_r(h)$, hence the limit  
\begin{equation}\label{eq10}
\lim_{h \to 0}{T(h)(f)}/{h^n}=K \text{ exists.}
\end{equation}
\\
If $n$ is even, then replacing $h$ with $-h$ in the above limit yields $-K=K$, or $K=0$.
If $n$ is odd, the condition $n>2m-1\geq 2k-1$ implies by Lemma~\ref{pr60} that $L=T_k(1)(x^n)\neq 0$. This in turn makes $T(1)(x^n)=T_k(2)(x^n)-2^{2k-1}T_k(1)(x^n)=(2^n-2^{2k-1})L=M\neq 0$. Consequently, the function
$$g(x)=\begin{cases}f(x)-(K/M)x^n & \text{  if $n$ is odd,}\\
f(x) & \text{  if $n$ is even,}
\end{cases}
$$
by \eqref{eq10}, has $T(h)(g)=o(h^n)$, regardless of the parity of $n$.
In order to apply Theorem~\ref{th2}, with $D(h)=T(h)(g)$, $R(h)=T_k(h)(g)$, and $s=2k-1$, we need to check that  $T_k(h)(g)=o(h^{2k-1})$. Indeed, this follows from the expression
\begin{equation}\label{eq6} 
T_k(h)(g)=\begin{cases}T_k(h)(f)-(KL/M)h^n  & \text{  if $n$ is odd,}\\
T_k(h)(f)& \text{  if $n$ is even,}
\end{cases}
\end{equation}
where the bound $2k-1\leq  2m-1\leq n-1$, together with the hypothesis $f(h)=o(h^{n-1})$, yields $f(h)=o(h^{2k-1})$, hence $T_k(h)(f)=o(h^{2k-1})$; we clearly have $h^n=o(h^{2k-1})$.
By Theorem~\ref{th2}, we deduce that $T_k(h)(g)=o(h^n)$ which, by~\eqref{eq6}, implies that the limit $$\lim_{h \to 0}\frac{T_k(h)(f)}{h^n}=\begin{cases}KL/M & \text{ if $n$ is odd,}\\0 & \text{ if $n$ is even,}\end{cases}$$ exists, as needed.
\end{proof}
\section{The Proof of the Generalized GGR Theorem}\label{S3}

As an application of the two versions of the second fundamental lemma, the following theorem provides sufficient conditions for an $n-1$ times Peano differentiable function at~$x$ to become $n$ times Peano differentiable at $x$.
This result is a generalization of the GGR theorem, since a consequence of it, Theorem~\ref{T8.4}, is a stronger version of the GGR theorem.

\begin{proposition}\label{cr1}
Let $n$ be an integer, at least $3$, and let $f$ be an $n-1$ times Peano differentiable function at $x$. If either
\begin{enumerate}
\item[(i)\,] $n$ is odd and  there are two $\mathcal{S}_{\frac{n+1}{2}}$-classes, one of odd and one of even differences of $f$ at $x$,~or
\item[(ii)] $n$ is even and there are an $\mathcal{S}_{\frac{n}{2}}$-class of even and an $\mathcal{S}_{\frac{n+2}{2}}$-class of odd differences of $f$ at $x$,
\end{enumerate}
then $f$ is $n$ times Peano differentiable at $x$. 
\end{proposition}
\begin{proof}
Let $g(h)=f(x+h)-p(x+h)$, where $p$ is the approximating Peano polynomial of $f$ at $x$. Then $g(h)=o(h^{n-1})$ and any difference $D_r$ in the above three classes has the property that $D_r(h)(g)$ evaluated at $0$ is the same as $D_r(h)(f-p)$ evaluated at $x$.

Notice that, for $n$ odd, the orders of even differences are at least $2(n+1)/2>n-1$, and the order of odd differences are at least $2(n+1)/2-1>n-1$; while for $n$ even, the orders of even differences are at least $2n/2>n-1$, and the orders of odd differences are at least $2(n+2)/2-1>n-1$. Hence all differences in Proposition~\ref{cr1} have orders greater than $n-1$, and so they vanish on $p$.
In particular, $D_r(h)(g)$ evaluated at 0 is the same as $D_r(h)(f-p)=D_r(h)(f)$ evaluated at $x$.

The conditions in both Theorems~\ref{th3} and~\ref{th4} are met for $g$, and so, by the same results, both limits $\lim_{h \to 0}[{g(h)+g(-h)}]/{h^n}$ and $\lim_{h \to 0}[{g(h)-g(-h)}]/{h^n}$ exist. Taking the average, $\lim_{h \to 0}{g(h)}/{h^n}$ exists, hence~$g$ is an~$n$ times Peano differentiable function at~$0$, or $f$ is $n$ times Peano differentiable at $x$.
\end{proof}
An interesting special case of the above theorem is $n=3$, when each of $\mathcal{S}_2^{\rm ev}$ and $\mathcal{S}_2^{\rm odd}$ consists of just one symmetric difference. Note that a third difference of a function~$f$ at~$x$ and~$h$ with nodes $-1,0,1,2$ is a non-zero scalar multiple of $\Delta_{3,1}(h)$, a third difference with nodes $-2,-1,0,1$ is a non-zero scalar multiple of $\Delta_{3,2}(h)$, and  a third difference with nodes $-2,-1,0,1,2$ is a non-zero linear combination of both, with the exception of $a\Delta_{3,1}(h)-a\Delta_{3,2}(h)=a\Delta_{4,2}(h)$, $a\neq 0$, or the fourth differences with these~nodes.

The following corollary shows that any third generalized Riemann derivative whose set of nodes is one of the above three is equivalent to the third Peano derivative, for all twice Peano differentiable functions $f$ at $x$.
\begin{corollary}
Let $D$ be any third difference of $f$ at $x$ and $h$ whose set of nodes is one of~{\small $\{-1,0,1,2\}$}, {\small $\{-2,0,1,2\}$}, or {\small $\{-2,-1,0,1,2\}$}.
If $f$ is twice Peano differentiable at $x$ and $\lim_{h \to 0}{D(h)(f)}/{h^3}$ exists,
then $f$ is three times Peano differentiable at $x$. 
\end{corollary} 
\begin{proof}
Take $\mathcal{S}_2^{\rm ev}=\{D(h)+D(-h)\}$ and $\mathcal{S}_2^{\rm odd}=\{D(h)-D(-h)\}$.  We need to check that none of the differences $D(h)\pm D(-h)$ is the zero difference, that is, $D(h)$ is not symmetric. Clearly, this is the case when the set of nodes is non-symmetric, $\{-1,0,1,2\}$ or $\{-2,0,1,2\}$. If the set of nodes is $\{-2,-1,0,1,2\}$, then 0 being a node makes the third difference $D(h)$ non-symmetric. The result now follows from Proposition~\ref{cr1}.
\end{proof}

The above corollary 1) proves both the classical and the updated versions of the GGR Theorem for $n=3$; and 2) shows that some crazy third generalized Riemann derivatives, such as
\[
-\tfrac 35\lim_{h\rightarrow 0}[2f(-2h)-7f(-h)+9f(0)-5f(h)+f(2h)]/h^3,
\]
obtained as $\tfrac 65\lim_{h\rightarrow 0}\Delta_{3,2}(h)/h^3-\tfrac 35\lim_{h\rightarrow 0}\Delta_{3,1}(h)/h^3$, are equivalent to the third Peano derivative for all twice Peano differentiable functions~$f$~at~0. These two important findings suggest that there should be a more general consequence of Proposition~\ref{cr1}, which will: 1) lead to a new proof of the GGR Theorem; and 2) provide large families of sets of $n$-th generalized Riemann derivatives that are equivalent to the $n$-th Peano derivative, for all Peano differentiable functions $f$ at $x$.

\medskip
This desired result for general $n$ is Theorem~\ref{T01}, the Generalized GGR Theorem, which we are ready to prove next.

\begin{proof}[Proof of the Generalized GGR Theorem]
Since $k\geq \lfloor n/2\rfloor $, two nodes of both differences $D _k(h)\pm D_k(-h)$ are $\pm k$ and the rest belong to the set $\{0,\pm 1,\dots,  \pm (k-1)\}$. Moreover, both limits $\lim_{h \to 0}[{D_k(h)(f)\pm D_k(-h)(f)}]/{h^n}$ exist.
We apply Proposition~\ref{cr1}, with
\[
\text{$\mathcal{S}_{\frac{n+1}{2}}^{\rm ev/odd}=\{D_k(h)\pm D_k(-h) : \tfrac{n+1}{2}\leq k \leq n-1\}$\;\qquad (for $n$ odd),}
\]
and with
\[
\quad\begin{aligned}
\mathcal{S}_{\frac{n}{2}}^{\rm ev}&=\{D_k(h)+ D_k(-h) : \tfrac{n}{2}\leq k \leq n-2\}\\
\mathcal{S}_{\frac{n+2}{2}}^{\rm odd}&=\{D_k(h)-D_k(-h) : \tfrac{n+2}{2}\,\leq \,k \,\leq \,n\}
\end{aligned}
\;\qquad\quad\text{(for $n$ even).}
\]
For this, we only need to verify the conditions in Definitions~\ref{d1} and~\ref{d2}, and then the result will follow from Proposition~\ref{cr1}.

First, none of $D_k(h)\pm D_k(-h)$ is the zero difference. Indeed, assuming the contrary, since $k$ is a node of $D_k(h)$, then $D_k(h)\pm D_k(-h)$ is zero implies that~$k$ is a node in~$ D_k(-h)$, or $-k$ is a node of $D_k(h)$. This could only happen when $k\leq \lfloor \frac n2\rfloor $, which for $n$ odd contradicts the condition $\frac {n+1}2\leq k$.
When $n$ is even, the same inequality will contradict the second assumption, $\frac {n+2}{2}\leq k$; and together with the first assumption, $\frac {n}{2}\leq k$, the same inequality forces $n=2k$,
so that  $D_{\frac{n}{2}}$ has order at least $n$, and its nodes belong to the set~$\{-\frac{n}{2},\ldots, \frac{n}{2}\}$ with $n+1$ elements.
This forces $D_{\frac{n}{2}}$ to be the even difference $\Delta_{n,\frac{n}{2}}$. Then $D_{\frac{n}{2}}(h)+D_{\frac{n}{2}}(-h)=2D_{\frac{n}{2}}(h)$ is non-zero, a contradiction.
 
Second, we need to check the orders condition in part (i) of Definitions~\ref{d1} and~\ref{d2}.
This means that for odd differences the order is at least $2\times \frac{n+1}{2}-1=n$ for $n$ odd and at least~$2\times \frac{n+2}{2}-1=n+1$ for $n$ even. For even differences we need to check that  the order is at least $2\times \frac{n+1}{2}=n+1$ for $n$ odd and at least $2\times \frac{n}{2}=n$ for $n$ even.
Summarizing, what we need to check is that, for $n$ odd, the order of $D_k(h)-D_k(-h)$ is at least $n$ and, for $n$ even, the order of $D_k(h)-D_k(-h)$ is at least $n+1$. The first of these is obvious, and the second follows from the fact that odd differences have odd orders. 

Finally, we need to check that, for $n$ odd, the order of $D_k(h)+D_k(-h)$ is at least $n+1$; and, for $n$ even, the order of $D_k(h)+D_k(-h)$ is at least $n$. The second of these is obvious, and the first follows from the fact that even differences have even order.
\end{proof}

In addition to providing infinitely many sets of $n$-th generalized Riemann derivatives whose joint existence is equivalent to the existence of the $n$-th Peano derivative, for all functions $f$ at $x$, the Generalized GGR Theorem provides a new proof of the GGR Theorem, which we restate here in its simplified, equivalent form, where half of the backward shifts of the $n$-th (forward) Riemann derivative are eliminated, as we discussed in the introduction. Recall that
$\Delta_{n,k}(h)(f)=\sum_{i=0}^n(-1)^i\binom nif(x+(n-k-i)h)$ is the $k$-th backward shift of the $n$-th Riemann difference $\Delta_{n}(h)(f)=\Delta_{n,0}(h)(f)$, and denote~$\mathscr{R}_{n,k}f(x)=\lim_{h\rightarrow 0}\Delta_{n,k}(h)(f)/h^n$. The notation $(n)$ means that the value $n$ is taken only for $n$ even.

\begin{theorem}[The Simplified GGR Theorem]\label{T8.4}
Let $n\geq 3$, and let $f$ be an $n-1$ times Peano differentiable function at $x$.

If all derivatives $\mathscr{R}_{n,k}f(x)$, for $k=\lfloor \frac{n+1}{2}\rfloor ,\lfloor\frac{n+1}{2}\rfloor+1,\ldots,n-1,(n)$, exist,
then~$f$ is~$n$ times Peano differentiable at~$x$.
\end{theorem}
\begin{proof}
The result follows from the Generalized GGR Theorem, by taking $D_k$ to be the difference $\Delta_{n,n-k}$, the $(n-k)$-th backward shift of $n$-th Riemann difference $\Delta_{n}$.
\end{proof}

\bibliographystyle{amsplain}

\begin{thebibliography}{AC}


\bibitem{As} J. M. Ash, \textit{Generalizations of the Riemann derivative,} Trans. Amer. Math. Soc. \textbf{126} (1967), 181--199. 


\bibitem{As2} J. M. Ash, \textit{Remarks on various generalized derivatives}. Special functions, partial differential equations, and harmonic analysis, 25--39, Springer Proc. Math. Stat., 108, Springer, Cham, 2014. 

\bibitem{AC} J. M. Ash and S. Catoiu, \textit{Quantum symmetric $L^p$ derivatives},
Trans. Amer. Math. Soc. \textbf{360} (2008), 959--987.

\bibitem{AC1} J. M. Ash and S. Catoiu, 
\textit{Multidimensional Riemann derivatives}, Studia Math. \textbf{235}  (2016),  no. 1, 87--100. 

\bibitem{AC2} J. M. Ash and S. Catoiu, \textit{Characterizing Peano and symmetric derivatives and the GGR conjecture's solution}, Int. Mat. Res. Notices IMRN 2022, no. 10, 7893--7921.

\bibitem{ACCh} J. M. Ash, S. Catoiu and W. Chin, \textit{The classification of generalized Riemann derivatives,} Proc. Amer. Math. Soc. \textbf{146} (2018), no. 9, 3847--3862.

\bibitem{ACCH} J. M. Ash, S. Catoiu and W. Chin, \textit{The classification of complex generalized Riemann derivatives,} J. Math. Anal. Appl. \textbf{502} (2021), no. 2, Article 125270. (40pp.) doi:10.1016/j.jmaa.2021.125270

\bibitem{ACCs} J. M. Ash, S. Catoiu and M. Cs\"{o}rnyei, \textit{Generalized vs. ordinary differentiation}, Proc. Amer. Math. Soc. \textbf{145} (2017), no. 4, 1553--1565. 





\bibitem{ACR}J.\ M. Ash, S. Catoiu, and R. R\'{\i}%
os-Collantes-de-Ter\'{a}n, \textit{On the nth quantum derivative,}\ J. Lond.
Math. Soc. \textbf{66} (2002), 114--130. 

\bibitem{AJJ} J. M. Ash, S. Janson, and R. L. Jones, \textit{Optical numerical differentiation using n function evaluations}, Calcolo \textbf{21} (1984), no. 2, 151--169.

\bibitem{AJ} J. M. Ash and R. L. Jones, \textit{Mean value theorems for generalized Riemann derivatives}, Proc. Amer. Math. Soc. \textbf{101} (1987), no. 2, 263--271. 




\bibitem{C} S. Catoiu, \textit{A differentiability criterion for continuous functions}, Monatsh. Math. \textbf{197} (2022), no. 2, 285--291.

\bibitem{dlVP} Ch. J. de la Vall\'ee Poussin, \textit{Sur l'approximation des fonctions d'une variable r\'eelle et de leurs d\'eriv\'ees par les p\^olynomes et les suites limit\'ees de Fourier}, Bull. Acad. Royale Belgique (1908), 193--254.

\bibitem{D}A. Denjoy, \textit{Sur l'int\'egration des coefficients diff\'erentiels 
d'ordre sup\'erieur,} Fund. Math. \textbf{25} (1935), 273--326.


\bibitem{EW} M. J. Evans and C. E. Weil, \textit{Peano derivatives: A survey}, Real Anal. Exchange \textbf{7} (1981-82), no. 1, 5--23. 

\bibitem{F} H. Fejzi\'c, \textit{Decomposition of Peano derivatives}, Proc. Amer. Math. Soc. \textbf{119} (1993), no. 2, 599--609.

\bibitem{F1} H. Fejzi\'c, \textit{Infinite approximate Peano derivatives}, Proc. Amer. Math. Soc. \textbf{131} (2003), no. 8, 2527--2536.

\bibitem{FFR} H. Fejzi\'c, C. Freiling, and D. Rinne, \textit{A mean value theorem for generalized Riemann derivatives}, Proc. Amer. Math. Soc. \textbf{136} (2008),  no. 2, 569--576.

\bibitem{FR} H. Fejzi\'c and D. Rinne, \textit{Peano path derivatives}, Proc. Amer. Math. Soc. \textbf{125} (1997), no. 9, 2651--2656.

\bibitem{FW1} H. Fejzi\'c and C. E. Weil, \textit{Repairing the proof of a classical differentiation result}, Real Anal. Exchange \textbf{19} (1993-94), 639--643.

\bibitem{FW} H. Fejzi\'c and C. E. Weil, \textit{A property of Peano derivatives in several variables}, Proc. Amer. Math. Soc. \textbf{141} (2013), no. 7, 2411--2417.

\bibitem{GGR} I. Ginchev, A. Guerraggio and  M. Rocca, \textit{Equivalence of Peano and Riemann derivatives.} Generalized convexity and optimization for economic and financial decisions (Verona, 1998), 169–178, Pitagora, Bologna,~1999. 

\bibitem{GGR1} I. Ginchev, A. Guerraggio and  M. Rocca, \textit{Equivalence of (n+1)-th order Peano and usual derivatives for n-convex functions}, Real Anal. Exchange \textbf{25} (1999/00),  no. 2, 513--520.

\bibitem{GR} I. Ginchev, M. Rocca, \textit{On Peano and Riemann derivatives}, Rend. Circ. Mat. Palermo (2) \textbf{49}  (2000),  no.~3, 463--480. 

\bibitem{HL} P. D. Humke and M. Laczkovich, \textit{Convexity Theorems for Generalized Riemann Derivatives}, Real Anal. Exchange \textbf{15} (1989/90), no. 2, 652--674.

\bibitem{HL1} P. D. Humke and M. Laczkovich, \textit{Monotonicity theorems for generalized Riemann derivatives}, Rend. Circ. Mat. Palermo (2) \textbf{38} (1989), no. 3, 437--454. 

\bibitem{Ki} A. Khintchine, \textit{Recherches sur la structure des fonctions mesurables,}
Fund. Math. \textbf{9} (1927), 212--279.

\bibitem{L} J. N. Lyness, \textit{Differentiation formulas for analytic functions}, Math. Comp. \textit{22} (1968), 352--362.

\bibitem{MZ} J. Marcinkiewicz and A. Zygmund, \textit{On the differentiability of functions and summability of trigonometric series,} Fund. Math. \textbf{26} (1936), 1--43.

\bibitem{MM} S. Mitra and S. N. Mukhopadhyay, \textit{Convexity conditions for generalized Riemann derivable functions}, Acta Math. Hungar.  \textbf{83}  (1999),  no. 4, 267--291.



\bibitem{P} G. Peano, \textit{Sulla formula di Taylor}, Atti Acad. Sci. Torino \textbf{27} (1891--92), 40--46.

\bibitem{RAA1} S. R\u{a}dulescu, P. Alexandrescu and D.-O. Alexandrescu, \textit{Generalized Riemann derivative},
Electron. J. Differential Equations 2013, No. 74, 19 pp. 

\bibitem{RAA2} S. R\u{a}dulescu, P. Alexandrescu and D.-O. Alexandrescu,
\textit{The role of Riemann generalized derivative in the study of qualitative properties of functions},
Electron. J. Differential Equations  2013, No. 187, 14 pp.

\bibitem{R} B. Riemann, \textit{\"Uber die Darstellbarkeit einer Funktion durch eine trigonometrische Riehe}, Ges. Werke, 2. Aufl., pp. 227--271. Leipzig, 1892.


\bibitem{Sa} H. E. Salzer, \textit{Optimal points for numerical differentiation}, Numer. Math. \textit{2} (1960), 214--227.

\bibitem{SZ} E. Stein and A. Zygmund, \textit{On the differentiability of functions}, Studia Math. \textit{23} (1964), 247--283.

\bibitem{T} B. S. Thomson, \textit{Monotonicity theorems}, Proc. Amer. Math. Soc. \textbf{83} (1981), 547--552. 


\bibitem{W} C. E. Weil, \textit{Monotonicity, convexity and symmetric derivatives}, Trans. Amer. Math. Soc. \textbf{231} (1976), 225--237.

\bibitem{Z} A. Zygmund, \textit{Trigonometric Series}, Vol. I, Cambridge University Press, 1959.
\end{thebibliography}

\end{document}